\def\a{\alpha}

\def\eps{\varepsilon}
\def\t{\theta}
\def\O{\Omega}

\def \dd#1{{\bf#1}}

\def\cl#1{{\cal#1}}

%Simbolos matematicos.

\def\id{\mathop{\rm id}\nolimits}

%Macros.

\def\ouv#1{\smash{\mathop{#1}\limits^{\lower 1pt\hbox
{$\scriptscriptstyle\circ$}}}}

\def\apl#1#2#3#4#5{{#1}:\matrix{
		\hfill#2 & \longrightarrow & #3\hfill \cr
		\hfill#4 & \longmapsto     & #5\hfill \cr}}

\def\hfl#1#2{\smash{\mathop{\hbox to 12mm{\rightarrowfill}}
\limits^{\scriptstyle#1}_{\scriptstyle#2}}}

%Titulos, enunciados.

\long\def\eno#1#2{\par\smallskip{\bf{#1}}{\it\ {#2}}\par\medskip}

\def\tit#1{\vskip 5mm plus 1mm minus 2mm {\tir #1}
		\vskip 3mm plus 1mm minus 2mm}

\def\stit#1{\vskip 3mm plus 1mm minus 2mm {\bf{#1}}
		\smallskip}

\long\def\dem#1#2{{\bf {#1}}{\ {#2}$\diamondsuit$}\medskip}

\font\tir=cmbx10 at 12pt

\def\ref#1#2#3#4{{\bf #1}{\ #2}{\it ,\ #3}{,\ #4}\medskip}

%Dibujos

\def \picture #1 by #2 (#3){\midinsert \centerline 
{\vbox to #2{\hrule width #1 heigth 0pt 
depth 0pt \null \vfill \special {picture #3}}}\endinsert }

\def\scaledpicture #1 by #2 (#3 scaled #4) {{
\dimen0 =#1 \dimen1 =$2
\divide \dimen0 by 1000 \multiply \dimen0 by #4
\divide \dimen1 by 1000 \multiply \dimen1 by #4
\picture \dimen0 by \dimen1 (#3 scaled $4)}}

\def\figure #1 #2 #3 {\midinsert \vglue 3mm 
{\vbox to #3 {\hrule width 6cm height 0cm depth 0cm \vfill
{\special {picture #1 scaled #2}}}}\vglue 2mm \endinsert}

\magnification=1200

%%%%%%%%%%%%%%%%%%%%%%%%%%%%%%
% Introduction
%%%%%%%%%%%%%%%%%%%%%%%%%%%%%%

\overfullrule=0pt

{\centerline {\tir {Linearization 
of holomorphic germs }}}

{\centerline {\tir {with resonant linear part}}}

\bigskip
\bigskip

{\centerline {{\bf R. P\'erez-Marco}\footnote {*} 
{  UCLA, 
Dept. of Mathematics, 
405, Hilgard Ave., Los Angeles, CA-90095-1555, 
USA, e-mail: ricardo@math.ucla.edu; CNRS UMR 8628, 
Universit\'e Paris-Sud, Math\'ematiques, 91405-Orsay, France.}}}

\bigskip

{\bf Abstract.} {\it We 
study the linearization problem of germs of 
holomorphic diffeomorphisms with resonant linear 
part. The formal linearization requires in general  
an infinite number of algebraic relations to be 
satisfied by the coefficients of the power series 
defining the holomorphic germ. For a fixed 
germ, the 
linearizing mappings are not unique, but the existence 
of a finite jet with a formal divergent linearization 
but no convergent one forces all  
the other linearizing mappings to be diverging.
We consider also polynomial families formally conjugated 
to a fixed resonant linear part. Then all the elements 
of the family are holomorphically linearizable, or we 
are in the precedent diverging situation for all but an 
exceptional set of parameter values of zero $\Gamma$-capacity.
In the case of maximal non-trivial resonance we prove the 
optimality of Bruno condition.
}

\bigskip

Mathematics Subject Classification 2000 :   37F50, 37J40, 32H50, 32U20, 41A17, 
34M35.

\medskip

Key Words : Small divisors, linearization, resonant germs, exceptional set. 

\bigskip
\bigskip

\tit {Introduction.}

The study of resonant germs of 
holomorphic diffeomorphisms $f: (\dd C^n , 0) \to (\dd C^n , 0)$,
$$
f(z)=A z +\cl O (z^2)
$$
is still in its infancy. 
The linear map $A\in GL_n (\dd C)$ is the differential of $f$ at 
$0$ and we can assume that it is given in Jordan normal form. 
We denote by $\{ \lambda_1 , \ldots , \lambda_n \}$ the set of 
eigenvalues of $A$ counted with multiplicity. The linear part 
$A$ is {\it resonant} if some relation of the form 
$$
\lambda_1^{i_1} \ldots \lambda_n^{i_n} -\lambda_j =0
$$ 
holds for some $i_1 ,\ldots , i_n \geq 0$, $|i_1 |+\ldots +|i_n|\geq 2$, 
and $1\leq j \leq n$.

A well studied resonant example (by H. Poincar\'e, G.D. Birkhoff, 
T.M. Cherry, C.L. Siegel, H. R\"ussman, H. Eliasson,...
[Po], [Bi], [Ch], [Si1], [Si2], [Si-Mo], [Ru], [El]) 
is the case of an elliptic fixed point at $0$
for a symplectic holomorphic map $f$.

In the non-resonant case, when $A$ is on the Poincar\'e domain, i.e.
$$
\max ( \max_i |\lambda_i | , \max_i |\lambda_i^{-1}| ) <1 \ ,
$$
it is well known that 
the holomorphic germ $f$ is holomorphically linearizable, that is,
there exists $h=\id +\cl O (z^2)$ a germ of holomorphic 
diffeomorphism, such that 
$$
h^{-1} \circ f \circ h =A \ .
$$

The Bruno condition on $A$ (or the eigenvalues of $A$), for a 
resonant or non-resonant $A$, is defined as
$$
-\sum_{k=1}^{+\infty } 2^{-k} \log  \Omega (2^k)  < +\infty
$$
where 
$$
\Omega (m)=\inf_{2\leq |i| \leq m \atop 1\leq j \leq n} |\lambda^i -\lambda_j| \ .
$$
and the infimum is taken among all $i$ such that 
$2\leq |i| \leq m$, $ \lambda^i -\lambda_j \not=0$, and $1\leq j \leq n$.

By the theorem of C.L. Siegel [Si3] and A.D. Bruno [Br], 
$f$ is also holomorphically linearizable  
when $A$ is non-resonant, diagonal, and the 
eigenvalues satisfy Bruno's diophantine condition stated above. 
 
\bigskip

When $f$ is formally linearizable (i.e. there exists as a formal map, 
defined by formal power series, conjugating formally $f$ to its linear part),
and the linear part satisfies Bruno's condition, then $f$ is also holomorphically 
linearizable (see [Br], [Ru], [El],  [E-V]). In this last situation the 
linearization is not unique.

When the formal linearization fails, one attempts to 
conjugate $f$ to a simpler normal form. In general the normalization mapping
is diverging (this has been well studied for the Birkhoff's normal 
form of a symplectic holomorphic germ with an elliptic fixed point).
In this paper, except for some common algebraic preliminaries, 
we study the linearization. The general conjugation to an arbitrary
normal form is treated in [PM2].

We consider the group of formal diffeomorphisms of $(\dd C^n ,0)$ 
(vanishing at $0$)
$$
\hat G=\hat {\hbox {\rm Diff}}(\dd C^n ,0) \subset \dd C^{\infty }
$$
We identify it to a subset of $\dd C^{\infty }$ by associating to 
$f\in \hat G$ the ordered set of the coefficients of its formal power 
series (choosing a monomial ordering). For $m\geq 1$, we consider 
the group $\hat G_m$ of $m$-jets of elements of $\hat G$, and 
the $m$-jet morphism
$$
\pi_m : \hat G\to \hat G_m \ .
$$
Let $d(m)$ be the number of coefficients of the $m$-jet.
We identify $\hat G_m$ to a subset of $\dd C^{d(m)}$. The group 
$\hat G_m$ is an affine algebraic group. We consider also the 
normal subgroup $\hat G_0 \subset \hat G$ of elements $h\in \hat G$
with identity linear part, $D_0 h=I$.

We have the following preliminary algebraic result, where we denote 
by $N$ a fixed normal form:

\eno {Theorem 1.}{Let $N\in \hat G$ and consider
the set $\hat G_N$ of elements $f\in \hat G$ such that 
$f$ is formally conjugated to $N$.
The subset $\pi_m (\hat G_N ) \subset \dd C^{d(m)}$ is an  
affine algebraic variety of dimension $d(N,m)$. Moreover, 
we have an onto polynomial map of degree less or equal than $m$
$$
\varphi_{m, N} : \dd C^{d(m)-m} \to \pi_m (\hat G_N ) 
$$ 
whose fibers are isomorphic non-singular affine algebraic varieties.

We assume $N=A\in {\hbox {\rm GL}}_n (\dd C )$.
For a given $f\in \hat G_A$, and any $m$, the $m$-jet projection of
the set $\hat G_0 (f)\subset \hat G_0$ 
of formal linearizations of $f$ is a linear sub-space 
of $\dd C^{d(m)-m}$ of dimension $\delta (A,m)=
d(m)-m-d(A,m)$. 
We have a linear isomorphisms
$$
\psi_{m,A,f} : \dd C^{\delta (A,m)} \to \pi_m (\hat G_0 (f)) \ .
$$
and a universal linearization linear isomorphism 
(obtained as inductive limit)
$$
\psi_{\infty , A, f} : \dd C^{\delta (A)} \to \hat G_0 (f) \ .
$$
where $\delta (A)=\lim_{m\to +\infty } \delta (A,m)$ (can be infinite).
}

We will need also an algebraic preliminary for polynomial 
families.

\eno {Theorem 2.}{ 
Let $(f_t)_{t\in \dd C^k}\subset \hat G_A$ be a polynomial 
family
$$
f_t (a)=Az +\sum_{0\leq |i|\leq d_0} t^i f_i (z)\ .
$$
We have a polynomial family of formal linearizations 
$(h_{t,s})_{(t,s)\in \dd C^k\times \dd C^{\delta (A)}}
\subset \hat G_I$ with 
$$
h_{t,s}= \psi_{\infty , A, f_t} (s) \ .
$$
More precisely, the coefficients of monomials of order $m$ 
of $h_{t,s}$ are polynomials on $t$ of degree less or equal to 
$d_0 m$ (and linear on $s$).}

After this algebraic preliminaries, we are ready to 
state the following theorems.
First for a single map.

\eno {Theorem 3.}{Let $f\in \hat G$ corresponding to a 
converging germ of holomorphic
diffeomorphism of $(\dd C^n , 0)$. Denote $G_0 (f)\subset \hat G_0(f)$ 
the set of converging linearizations of $f$.

We have the following dichotomy:

\item {1)} Up to finite order, all formal linearizations in 
$\hat G_0 (f)$ are converging. More precisely, for any $m\geq 2$, we have
$$
\pi_m (\hat G_0 (f))=\pi_m (G_0 (f)) \approx \dd C^{\delta (A,m)}\ ,
$$
and the set of converging linearizations is in one-to-one 
correspondence with the holomophic centralizer of $A$.
Morever, there is a uniform radius of convergence $R_0 >0$
such that for any $m\geq 2$ and $s_m \in \dd C^{\delta (A,m)}$
there exists a linearization $h$ of $f$ with $\pi_m (h) =s_m$
and radius of convergence at least $R_0$.

\item {2)} All formal linearizations in $G_0 (f)$ are 
diverging. 
}

After this theorem, we can talk without ambiguity about $f$ 
being linearizable or not. We say that $f$ is linearizable if
it has a converging linearization, thus we are in case (1), 
and this notion is essentially independent of the 
linearization chosen.

J.-Ch. Yoccoz in [Yo] appendix 1 studies the 
linearization of germs with non-diagonalizable
linear part. He proves that in this situation there 
always exists non-linearizable examples. More precisely,
that in a one dimensional family with linear dependence
on the parameter for almost all values of the parameter 
the dynamics is not linearizable. Yoccoz's argument 
is incomplete in the resonant case. Next theorem (see also [PM1])
shows that the set of exceptional parameter values indeed
has capacity $0$ and the result includes the resonant case. 

\medskip

The next theorem is guided by the general principle stated in 
[PM1] that in generic families of dynamical systems
presenting problems of small divisors :

\bigskip

{\bf We have total convergence or general divergence except for a very small exceptional set in parameter space}

\bigskip

We use below the notion of $\Gamma$-capacity. We recall its definition
in section 2 and 
refer to [Ro] for basic properties.

\eno {Theorem 4.}{Let $(f_t)_{t\in \dd C^k}$ be a polynomial family as
above, and let $(h_{t,s})_{(t,s)\in \dd C^k \times \dd C^{\delta (A)} }$
be the universal family of uniformizations. We have the following 
dichotomy:

\item {1)} Up to finite order, all the formal linearizations 
$(h_{t,s})$ are converging, i.e. all $f_t$ are linearizable.
More precisely, for any $m\geq 2$ and for any $(t, s_m)\in \dd C^k \times 
\dd C^{\delta (A,m)}$, there exists $s\in \dd C^{\delta (A)}$, with
$s_m=\pi_m (s)$ and $h_{t,s}$ converging.

\item {2)} Except for an exceptional set $E$ in $\dd C^k$ of 
$\Gamma$-capacity $0$, all formal linearizations in $h_{t,s}$, $t
\notin E$, are diverging, i.e. all $f_t$, $t\notin E$, are non-linearizable.
}

\bigskip

In the last part of the article we study some examples. 
We give examples
of resonant linear part $A$ with $\delta (A)=+\infty$.

\medskip

A measure of the importance of resonance is the 
$\dd Q$-torsionless dimension of the $\dd Q $-vector space
$$
V=V(\lambda_1 , \ldots \lambda_n)= \dd Q \oplus 
\bigoplus_{j=1}^n  \dd Q  \ \left (
{1\over 2\pi i} \log \lambda_j \right )   
$$ 
that is
$$
l(A)=\dim_{\dd Q} V /\dd Q  \ .
$$

The linearization problem in the case of minimal dimension 
$l(A)=0$ and $A$ is diagonal is simple.  
Linearization is equivalent
to finite order (for the composition) of 
the germ. 
The next non-trivial case corresponds
to $l(A)=1$. Then we have two cases, depending on wether
$$
\dd Q /\dd Z \cap 
\bigoplus_{j=1}^n (\dd Q /\dd Z ) \ 
\left ( {1\over 2\pi i} \log \lambda_j \right )  =\{ 0 \}
$$
(torsionless) or the intersection is not reduced to $0$ (torsion).

There exists $\mu=e^{2\pi i \t }$ and a 
$q=q(A)$-root of unity $\eps$ ($\eps =1$ in the torsionless case,
$q\geq 2$ in the torsion case) 
such that 
$\lambda_i =\eps^{a_i} \mu^{b_i}$ for some 
integers $a_i$ and $b_i$. The following result seems
to have been unnoticed (but see [Go] in for symplectic holomorphic
maps in dimension $2$):

\eno {Theorem 5.}{ The eigenvalues $(\lambda_1 , \ldots , \lambda_n )$
do satisfie Bruno condition if and only if $\mu$ satisfies the 
one dimensional Bruno condition (i.e. $\t $ is a Bruno number).

The optimal arithmetic condition of 
holomorphic linearization of resonant formally 
linearizable holomorphic germs 
with linear part $A$ with $l(A)=1$ is Bruno condition.}

At this respect one should note that M. Herman 
has conjectured ([He] conjecture 1 p.147)
in the non-resonant and diagonal case, that Bruno condition is optimal in all 
dimensions. Yoccoz has verified this in dimension $1$ ([Yo]). The author
conjectures that Herman's conjecture is false in dimension larger 
than $2$.

\medskip

The torsion case is irrelevant for purposes of linearization.
Observe that 
$$
\dd Q \cap \bigoplus_{j=1}^n 
\dd Z \ \left ( {1\over 2\pi i} \log \lambda_j \right ) 
$$
is a $\dd Z$-submodule of a free finitely generated 
$\dd Z$-module, thus it is finitely
generated since $\dd Z$ is Noetherian.
If there is torsion, the torsion 
$q=q(A)$ of $(\lambda_1 , \ldots , \lambda_n)$ is the 
minimal positive integer $q\geq 2$ (which exists  
by the precedent observation) such that 
$$
{1\over q} \dd Z \supset \dd Q \cap \bigoplus_{j=1}^n 
\dd Z \ \left ( {1\over 2\pi i} \log \lambda_j \right ) \ .
$$
If $f$ is a resonant germ with linear part
$A$ with torsion $q\geq 2$ then $f^q$  has linear part 
$A^q$ with no torsion. Moreover, it is easy to see that 
$f$ is linearizable if and only if $f^q$ is linearizable.

The case $l(A) \geq 2$ deserves further investigations
that will be pursued in the future.

\bigskip

\null
\vfill
\eject

%%%%%%%%%%%%%%%%%%%%%%%%%%%%%
%  Proof
%%%%%%%%%%%%%%%%%%%%%%%%%%%%

\stit {1) Algebraic preliminaries.}

The study of conjugacy to a non-linear normal form
is done in [PM2]. Part of the formal theory is the 
same as for the problem of linearization, this is why we 
include a general first part in theorem 1.
We don't develop here a complete algebraic study.
For basic notions relative to algebraic groups we refer 
to [Sh]. Most of the work done here can be 
straithforward generalized 
to more general fields of coefficients than $\dd C$.

Using the notations of the introduction, we consider
a normal form $N\in \hat G$ with   
$$
D_0 N =A
$$
and its (formal) centralizer in $\hat G_0$,
$$
{\hbox {\rm Cent }} (N) =\{ h\in \hat G_0 ; h\circ N =N\circ h \}\ .
$$
Note that for each $m\geq 2$, 
$$
\pi_m ({\hbox {\rm Cent }} (N) ) \subset \pi_m (\hat G_0 )
\subset  \dd C^{d(m)}
$$
is an affine algebraic subgroup. Thus it can be identified 
to a linear group.

It is  natural to ask
\bigskip
  
{\bf Question :}
{\it Which algebraic groups do arise as
groups ${\hbox {\rm Cent }} (N)$ for some $N$ ?
}
\bigskip

This question does not seem to have been explored at all.

Given $f\in \hat G$, we consider the set  
$$
\hat G_0 (f,N) =\{ h\in \hat G_0 ; h^{-1} \circ f\circ h =N \} \ .
$$
of normalizations of $f$.

\eno {Proposition 1.1.}{For each $m\geq 2$,
$$
\pi_m (\hat G_0 (f,N) )\subset \dd C^{d(m)}\ 
$$
is a smooth affine algebraic variety.
}

\dem {Proof.}{ If the set $\hat G_0 (f,N)$ 
is empty there is nothing to 
prove. Otherwise choose $h_0 \in \hat G_0 (f,N)$. The map
$$
\apl {\Phi } {\pi_m ({\hbox {\rm Cent }} (N) )}
{\pi_m (G_0 (f,N))} {\pi_m (h)} {\pi_m (h\circ h_0)}
$$
is a linear isomorphism. The inverse is simply given
by
$$
\pi_m (h)\mapsto \pi_m (h\circ h_0^{-1} ) \ .
$$
Now $\pi_m ({\hbox {\rm Cent }} (N))$ is an affine algebraic 
group, thus it is non-singular.}

We consider now the set $\hat G_N$ of elements $f\in \hat G$
which are formally conjugated to $N$,
$$
\hat G_N =\{ h^{-1}\circ N\circ h ; h\in \hat G_0 \}
$$
If $N$ has a resonant linear part $A$ 
we have $\hat G_N \not= \{f\in \hat G ; D_0 f=A \}$. 
Fixing $m\geq 2$, 
we can identify $\pi_m (\hat G_0 )$ to $\dd C^{d(m)-m}$.
The map 
$$
\varphi_{m,N} : \dd C^{d(m)-m} \to \pi_m (\hat G_N )
$$
defined by $\pi_m (h) \mapsto \pi_m (h^{-1} \circ N \circ h)$
is a polynomial map because the coefficients of the power 
series defining $h^{-1}$ depend polynomially on those of 
$h$. It is not difficult to check that each 
coordinate function is a polynomial
in $d(m)-m$ variables of degree less than $m$.
Each of its fibers is isomorphic as
algebraic varieties to $\pi_m ({\hbox {\rm Cent }} (N))$.
Thus it is a non-singular affine variety.
To end the proof of theorem 1 we still have to construct 
the map $\psi_{m,A,f}$. We do this in the case where 
$N=A$ is linear. 

\eno {Lemma 1.2.}{Let $m\geq 2$.
The $m$-jet of the centralizer of $A$, 
$\pi_m ({\hbox {\rm Cent }}(A)) \subset \pi_m (\hat G_0)
\subset \dd C^{d(m)-m}$, is a linear sub-space of 
dimension $\delta (A,m)\leq d(m)-m$.}

\dem {Proof.}{Let $h\in \hat G_0 \subset 
\dd C^{d(m)-m}$. Identifying
the coefficients of monomials in the equation 
$$
\pi_m (h\circ A)=\pi_m (A\circ h)
$$
we get $d(m)-m$ linear equations in the 
$d(m)-m$ coefficients of $\pi_m (h)$.}

The maps $\psi_{m,A}: \dd C^{\delta (A,m)}\to 
\pi_m ({\hbox {\rm Cent (A)}})$ are defined choosing 
a bases for ${\hbox {\rm Cent (A)}}$.
Finally the map $\psi_{m,A,f}$ is defined by 
$$
\psi_{m,A,f} =\Phi\circ \psi_{m,A}
$$
where $\Phi$ is the linear map used in the proof of 
proposition 1.1. Note that the construction is not 
canonical since it depends on the choice of the 
bases of ${\hbox {\rm Cent (A)}}$ and a choice of 
$h_0\in \hat G_0(f,A)=\hat G_A$ for $\Phi$.

\bigskip

Theorem 2 for a polynomial family of formal germs $(f_t)$ 
follows from an elimination argument. Writting down the 
equations for the linearization of $f_t$, at each order $m$
the equations determine some coefficients of order $\leq m$ 
of the linearization in terms of $t$. We set the other coefficients
equal to $0$. The system of equations are compatible by assumption.
In that way we obtain a polynomial family $(h_t)$ of linearizations.
One easily checks that the coefficient of a monomial of order
$m$ is a polynomial on $t$ of degree at most $m d_0$. Now for each 
$t$ we have a map $\Phi_t : {\hbox {\rm Cent}}(A) 
\to \hat G_0(f_t , A)$, $h\mapsto h\circ h_t$. We define
$$
\psi_{\infty , A, f_t}=\Phi_t \circ \psi_{\infty , A}
$$
and $h_{t,s}=\psi_{\infty , A, f_t}$.

\medskip

One can do the above constructions in a more "functorial"
way, but this is useless for our purposes.

\stit {2) Proof of theorem 2 and 3.}

\stit {Proof of theorem 2.} 

The proof follows from theorem 1. Assume that we have 
a converging linearization $h_0$. From theorem 1 we know 
that all other linearizations are of the form $h\circ h_0$ 
where $h$ runs over all elements of ${\hbox {\rm Cent}}(A)$.
In particular, choosing $h$ to be polynomial of order at 
most $m$ we see that 
$$
\pi_m (\hat G_0 (f))=\pi_m (G_0 (f)) \ .
$$

\stit {Proof of theorem 3.}

We recall Berstein's lemma in approximation theory 
(see [Ra] p.156): 

\eno {Lemma (Bernstein).}{Let $K\subset \dd C$ be a non-polar 
set, and $\O$ be the component of ${\overline {\dd C}}-K$ containing
$\infty$.

If $P$ is a polynomial of degree $n$, then for $z\in \dd C$
$$
|P(z)|\leq ||P||_{C^0(K)} e^{ng_{\O} (z, \infty )}
$$
where $g_{\O }$ is the Green function of $\O$.}

\stit {$\Gamma$-capacity.}

We recall the definition of $\Gamma$-capacity and 
we refer to [Ro] for more properties.
Let $E \subset \dd C^k$. The $\Gamma$-projection of $E$ on 
$\dd C^{k-1}$ is the set $\Gamma_k^{k-1} (E)$ 
of $z=(z_1, \ldots , z_{k-1})\in \dd C^{k-1}$
such that 
$$
E\cap \{ (z, w)\in \dd C^k \}
$$
has positive capacity in the complex plane 
$ \dd C_z=\{ (z, w)\in \dd C^k \}$.
We define 
$$
\Gamma_k^1 (E)=\Gamma_2^1 \circ \Gamma_3^2\circ \ldots \Gamma_k^{k-1} (E) \ .
$$
Finally, the $\Gamma$-capacity is defined as 
$$
\Gamma {\hbox {\rm -Cap}} (E)=\sup_{A\in U(k,\dd C)} {\hbox {\rm Cap}}
\ \Gamma_k^1 (A(E)) \ .
$$
where $A$ runs over all unitary transformations of $\dd C^k$.

Using the definition of $\Gamma$-capacity it is easy to see 
that we are reduced to prove theorem 3 for $k=1$

\bigskip

From the algebraic preliminaries it follows 

\eno {Lemma 2.1.}{The coefficient vectors $h_i(t, s_i)$, 
$s_i =\pi_i (s)$ of the
formal linearization  
$$
h_{t, s} (z)=z+\sum_{i=(i_1, \ldots i_n)\atop 
i_1+\ldots +i_n \geq 2} h_i(t, s_i) z^i
$$
have coordinates that are linear in the coordinates 
of $s_i$ and polynomial in the parameter 
$t=(t_1, \ldots t_k)$ 
of degree less than $d_0 (i_1+\ldots i_n)$.}

Taking into account the definition of $\Gamma$-capacity
given above, as observed before,  we can assume that  
$k=1$.

Let 
$$
E=\{ t\in \dd C ; f_t \ {\hbox {\rm    is linearizable  }} \} \ .
$$
Fix $m\geq 2$, $s_m \in \dd C^{\delta (A , m)}$ and
$s=(s_m , 0,\ldots )$.
We want to show that $E$ is polar or the whole complex plane.
We have 
$$
E=\bigcup_{j\geq 1} E_j
$$
where $E_j$ the set of parameters $t$ such that $h_{t,s}$ has radius
of convergence larger or equal to $1/j$. Thus if $E$ is non-polar, 
we have that for some 
$j\geq 1$, $E_j$ is not polar. Thus there exists $\rho_0 >0$
such that for all $t\in E_j$, 
$$
\varphi (t)=\limsup_{|i|\to +\infty} ||h_i(t,s)||\rho_0^{-|i|} <+\infty \ .
$$
The function $\varphi$ is lower semicontinuous, and
$$
E_j=\bigcup_p L_p
$$
where $L_p=\{ z\in E_j ; \varphi (t) \leq p\} $ is closed.  
By Baire theorem for some $p$, $L_p$ has non-empty interior (with 
respect to $E_j$), thus this $L_p$ has positive capacity. Finally
we found a compact set $C=L_p$ of positive capacity such that 
there exists $\rho_1 >0$ such that for any $t\in C$ and 
and all $i\in \dd N^n$,
$$
||h_i (t,s)||_{C^0 (C)} \leq \rho_1^{|i|} \ .
$$
Using Bernstein's lemma and lemma 1.1 we get that for any compact 
set $K\subset \dd C$ we have
$$
||h_i(t)||_{C^0 (K)} \leq C(K)^{d |i|} \rho_1^{|i|} \ ,
$$
for some constant $C(K)$ depending only on $K$.
Thus $f_t$ is linearizable for any $t\in \dd C$.
The constant $C(K)$ can be estimated by the precise form 
of Bernstein lemma as
$$
C(K) = \exp (\sup_K g_{\O} (t , \infty ) )
$$
where $\O$ is the connected component containing $\infty$ of
the complement of $C$. The asymptotic
$$
g_\O (t, \infty) \approx \log |t|
$$
for $t\to \infty$ can be used to give a lower estimate on the 
radius of convergence.

\stit {3) Some examples.}

\stit {a) Examples with $\delta (A)=+\infty$.}

{\bf Elliptic point.}

A well known example is an elliptic 
fixed point. For example, consider
$$
A=\pmatrix {\lambda &0 \cr 0 &\lambda^{-1} \cr }
$$
with $\lambda \in \dd C^*$ and $\lambda$ not a 
root of unity.

\eno {Lemma 3.1.}{ The The holomorphic (resp. formal) 
centralizer in $G_0$ (resp. $\hat G_0$) of $A$ is composed 
by maps of the form 
$$
(z_1 ,z_2) \mapsto l(z_1 ,z_2)=(z_1 +z_1 \varphi_1 (z_1z_2) , 
z_2+z_2 \varphi_2 (z_1 z_2))
$$
where $\varphi_j(z)=\sum_{i=1}^{+\infty } \varphi_{j, i} z^i$ is a 
holomorphic (resp. formal) 
power series.}

\dem {Proof.}{The proof is straightforward identifying 
coefficients in $l\circ A =A\circ l$.}

\stit {Non-trivial Jordan block.}

We consider now another more elaborate example.
Let 
$$
A=\pmatrix {1 & 1 \cr 0 & 1 \cr }
$$

\eno {Proposition 3.2.}{The holomorphic (resp. formal) 
centralizer in $G_0$ (resp. $\hat G_0$) of $A$ is composed 
by maps of the form 
$$
(z_1 ,z_2) \mapsto l(z_1 ,z_2)=(z_1 +k(z_2) , z_2)
$$
where $k(z)=\sum_{i=2}^{+\infty } k_i z^i$ is a 
holomorphic (resp. formal) 
power series.}

\eno {Lemma 3.3.}{Let $\psi (z_2)$ 
be a given formal 
power series of valuation $\geq 2$. We consider the 
equation 
$$
\varphi (z_1+z_2 , z_2) =\varphi (z_1 , z_2) +\psi (z_2)\leqno {(*)}
$$
and seek solutions $\varphi (z_1 , z_2)$ which are 
formal power series of valuation $\geq 2$. 

If $\psi \equiv 0$ then the solutions of $(*)$ are 
the formal power series $\varphi (z_1 , z_2)= \rho (z_2)$
independent of $z_1$.

If $\psi $ is not identically $0$ there are no solutions 
to equation $(*)$.}

\dem {Proof.}{We consider the linear operator
$$
L: \varphi (z_1 , z_2) \mapsto \varphi (z_1+z_2 , z_2)-
\varphi (z_1 , z_2)
$$
defined in the vector space $E$ of formal power series in the
two variables $(z_1 ,z_2)$ with valuation $\geq 2$. This 
linear operator leaves invariant the finite dimensional vector
space $E_n$, $n\geq 2$, of homogeneous polynomials of degree $n$.
In order to solve $(*)$ we are reduced to solve the equation in $E_n$
for $n\geq 2$. 
Choosing the base $(z_1^n , z_1^{n-1} z_2 , \ldots , z_2^n )$
of $E_n$, the matrix of $L_n=L_{/E_n }$ is triangular with combinatorial 
coefficients $\pmatrix {i \cr j\cr }$, 
$$
M=\pmatrix { 0 & 0& 0 & \ldots   & 0 \cr
1 & 0 & 0 & \ldots & 0 \cr
1 & 2 & 0 & \ldots & 0\cr
\vdots & \vdots & \vdots & \ddots& \vdots  \cr
1 & n & \ldots  & n& 0\cr
}
$$
The only eigenvalue is $0$, the kernel is
spanned by  
$(0, \ldots , 0, 1)$, and the image  
by the first $n-1$ vectors of the bases. 
So  
the $n$-homogeneous part of $\psi (z_2)$ is not zero, 
then it doesn't belong to the image of $L$, thus $(*)$ 
has no solution. Now if $\psi (z_2)\equiv 0$ then the 
homogeneous solutions of $(*)$ is the kernel
of $L_n$.
So a solution $\varphi $ of $(*)$ only depends on $z_2$.}

\dem {Proof of proposition 3.2.}{Obviously any map of the
above form commutes 
with $A$. Conversely, we write $l(z)=(z_1 +l_1(z) , z_2 +l_2 (z))$
with $l_i$ a formal power series in $(z_1, z_2)$ of 
valuation $\geq 2$.
The second coordinate in the equation 
$$
l\circ A = A \circ l
$$
gives (eliminating linear parts)
$$
l_2(z_1+z_2 , z_2) =l_2 (z_1, z_2) \ .
$$
This implies by the previous lemma that $l_2 (z_1 , z_2)= \psi (z_2)$.
Now we look at the first coordinate of the commuting equation.
We get
$$
l_1 (z_1 +z_2 , z_2) =l_1 (z_1 , z_2) +l_2 (z_1 , z_2)=
l_1 (z_1 , z_1)+\psi (z_2)\ .
$$
Using the lemma again we get that $\psi$ must be zero and 
$l_1 (z_1 , z_2) = k(z_2)$.q.e.d.
}

Along the same lines one can determine the centralizer of
a non-trivial Jordan block
$$
A=\pmatrix {1 & 1 & 0 & 0 & \ldots & \ldots & \ldots  & 0 \cr
0&1&1&0&\ldots &\ldots &\ldots &\vdots \cr
\vdots & \ddots & \ddots & \ddots & \ddots &\ldots &\ldots &\vdots \cr
\vdots &\vdots &\ddots & 1& 1& 0 & \ldots &\vdots \cr
\vdots &\vdots &\vdots &\ddots & 1 &0 &\ddots & \vdots \cr
\vdots &\vdots &\vdots &\vdots &\ddots &\ddots &\ddots & \vdots \cr
\vdots &\vdots &\vdots &\vdots &\vdots &\ddots &1 &0 \cr
0& \ldots  & \ldots & \ldots & \ldots  & \ldots  & 0 & 1 \cr
}
$$
which is also infinite. 

\eno {Proposition 3.4}{If the matrix $A$ contains a 
non-trivial Jordan block with an eigenvalue $1$ (or a root of 
unity) then $\delta (A) =+\infty$.}

\stit {b) Linearization of resonant germs with $l(A)=1$.}

All eigenvalues $\lambda_j$ are non zero. We write 
$\lambda_j =e^{2\pi i \a_j}$ for some $\a_j \in \dd C$
determined up to an additive integer. We write
$$
(\a_1 , \ldots , \a_n )\in \cl B
$$
to express that $(\lambda_1 , \ldots , \lambda_n)$ 
satisfies the Bruno condition.

With the notations of the introduction, we assume that $l(A)=1$,
that is
$$
V/\dd Q =\dd Q \ \t \ .
$$
Then there is a root of unity $\eps$ and integers $(a_k)$
and $(b_{k})$ such that 
$$
\lambda_k =\eps^{a_k}  \mu^{b_{k}}
$$
where $\mu =e^{2\pi i \t }$.
We prove

\eno {Proposition 3.5.}{We have 
$$
(\a_1 , \ldots , \a_n )\in \cl B
$$
if and only if $\t \in \cl B $.
}

\eno {Lemma 3.6.}{Given a positive integer $q\geq 1$, 
we have 
$$
(\a_1 , \ldots , \a_n )\in \cl B
$$
if and only if 
$$
(q\a_1 , \ldots , q\a_n )\in \cl B \ .
$$
}

\dem {Proof.}{We observe that if $|a|=|b|=1$ and 
$a^q-b^q$ is small then 
$$
|a^q-b^q|\approx q|a-b| \ .
$$
Using this, we have that the small divisors for $\a$ and 
$q \a$ are the same, more precisely, for $m\geq 2$,
$$
\O_{q\a} (m)\approx \O_{\a }(m)
$$
and the result follows.}

\dem {Proof of proposition 3.5.}{Using the lemma 
we are reduced to prove the theorem in the case 
without torsion, i.e. $\eps=1$. 

Observe that a small divisor for $(\lambda_1, \ldots , \lambda_n)$
can be written as
$$
\lambda_1^{i_1}\ldots \lambda_n^{i_n} -\lambda_j
=\mu^{\sum_k b_k i_k}-\mu^{b_j}=\mu^{b_j} \left (
\mu^{\sum_k b_k i_k -b_j} -1 \right )
$$
Also $|b.i|\leq ||b||_{\infty } |i|$, thus
$$
\O_{\a } (m) \geq C \ \O_{\t } (||b||_{\infty } m) \ .
$$
Also if we have a small divisor for $\mu$, $\mu^n-\mu$, 
and $b_j\not=0$ then 
$$
\lambda_j^n-\lambda_j=\mu^{nb_j }-\mu^{b_j }\approx b_j (\mu^n -\mu )
$$
so we conclude that 
$$
\O_{\t } (m)\geq C \ \O_{\a } (m) \ .
$$
}

\stit {Proof of theorem 5.}

The first assertion has been proved above. 
The linearization under Bruno condition follows 
from the general theory. When Bruno condition is 
violated, then according to lemma 3.3 we have that 
all non-rational $\a_j$ do not satisfie Bruno condition. We just pick 
one such $\a_j$ and a non-linearizable holomorphic germ in one variable
$$
g(z)=e^{2\pi i \a_j } z +\varphi (z)
$$
with $\varphi (z)=\cl O (z^2 )$.
We construct $f:(\dd C^n ,0) \to (\dd C^n , 0)$ by
$$
f(z_1, \ldots  , z_n )=Az+(0, \ldots , 0, \varphi (z_j) , 0 , \ldots, 0)
$$
where we can assume that the $z_j$-axes is an eigendirection for $A$.
Then $f$ is formally linearizable but not analytically linearizable.

\bigskip

We finish the proof of the last comments in the introduction
about the irrelevance of the torsion for the purposes of 
linearization.

\eno {Proposition 3.7.}{Let $q\geq 1$. The germ $f$ is 
holomorphically (resp. formally) linearizable if and only
if $f^q$ is holomorphically (resp. formally) linearizable.}

\dem {Proof.}{We prove the non-trivial statement.
If $f^q$ is holomorphically (resp. formally) 
linearizable there an element
$k\in G_0$ (resp. $k\in \hat G_0$) such that 
$$
k\circ f^q=A^q\circ k \ .
$$ 
Then 
$$
k_0 ={1\over q} \sum_{i=0}^{q-1} A^{-i} \circ k\circ f
$$
is an element of $G_I$ (resp. $\hat G_I$) that linearizes 
$f$.}

\eno {Proposition 3.8.}{If $A=D_0 f$, $l(A)=0$ and 
$A$ is diagonal then 
$f$ is linearizable if and only if 
$$
f^{q(A)} =\id \ .
$$
}

\dem {Proof.}{We have $A^{q(A)}=I$ so if $f$ is linearizable
we get $f^{q(A)}=\id $.  The converse results from the 
application of the previous proposition.}

%%%%%%%%%%%%%%%%%%%%%%%%%%%

\null
\vfill
\eject

\bigskip

{\centerline {\bf BIBLIOGRAPHY}}

\bigskip
\bigskip

\ref{[Br]}{A.D. BRJUNO}{Analytical forms of differential equations}
{Trans. Mosc. Math. Soc., {\bf 25}, 1971, {\bf 26}, 1972.}

\ref{[Ch]}{T.M. CHERRY}{On the solution of hamiltonian systems of 
differential equations in the neighborhood of a singular point}
{.}

\ref{[El]}{H. ELIASSON}{Hamiltonian systems with linear 
form near an invariant torus}
{Non-linear Dynamics, Bologna, 1988, World Sci. Publishing,
Teaneck, NJ, 1989, p.11-29.}

\ref{[E-V]}{J. ECALLE, B. VALLET}{Correction and linearization 
of resonant vector fields and diffeomorphisms}
{Math. Z., {\bf 229}, 1998, p.249-318.}

\ref{[Go]}{X. GONG}{Conformal maps and non-reversibility of elliptic 
area preserving maps}
{to appear in Math. Z.}

\ref{[He]}{M. R. HERMAN}{Recent results and some 
open questions on Siegel's lin\'earisation theorem 
of germs of complex analytic diffeomorphisms 
of $\dd C^n$ near a fixed point}{Proceedings 
$VIII^{th}$ Int. Conf. Math. Phys., World Scientific 
Publishers, Singapour, 1987.}

\ref{[PM1]}{R. P\'EREZ MARCO}{Total convergence or general divergence
in Small Divisors} 
{Preprint, 2000.}

\ref{[PM2]}{R. P\'EREZ MARCO}{On the convergence of 
Birkhoff's normal form}{Preprint 2000.}

\ref{[Po]}{H. POINCAR\'E}{Les m\'ethodes nouvelles de la m\'ecanique
c\'eleste}{volume 1, Paris, 1893.}

\ref{[Ra]}{T. RANSFORD}{Potential theory in the complex plane}{London
Mathematical Society, Student Texts {\bf 28}, Cambridge University 
Press, 1995.}

\ref{[Ro]}{L.I. RONKIN}{Introduction to the theory of 
entire functions of several variables}{Translations of Mathematical
monographs, American Mathematical Society, {\bf 44}, 1974.}

\ref{[Ru]}{H. R\"USSMANN}{On the convergence of power series 
transformations of analytic mappings near a fixed point}
{Preprint I.H.E.S., Paris, 1977.}

\ref{[Sh]}{I.R. SHAFAREVITCH}{Basic algebraic geometry}{Springer Verlag, 1994.}

\ref{[Si1]}{C.L. SIEGEL}{On the integrals of canonical systems}{Annals
of Mathematics, {\bf 42}, {\bf 3}, 1941, p.806-822.}

\ref{[Si2]}{C.L. SIEGEL}{\"Uber die Existenz einer Normalform
analytischer Hamiltonscher Differentialgleichungen in der 
N\"ahe einer Gleichgewichtl\"osung}{Math. Annalen, {\bf 128},
1954, p.144-170.}

\ref{[Si3]}{C.L. SIEGEL}{Iterations of analytic 
functions}{Ann. Math.,{\bf 43}, 1942, p. 807-812.}

\ref{[Si-Mo]}{C.L. SIEGEL, J. MOSER}{Lectures on celestial mechanics}
{Springer-Verlag, {\bf 187}, 1971.}

\ref{[Yo]}{J.-C. YOCCOZ}{Th\'eor\`eme de Siegel, 
polyn\^omes quadratiques et nombres de Brjuno}
{Ast\'erisque {\bf 231}, 1995.}

\end